\documentclass[a4paper]{article}
\usepackage{amssymb}
\usepackage{amsmath}
\usepackage[authoryear]{natbib}
\usepackage[pdftex]{graphicx,color}
\usepackage{authblk}

\newcommand{\bs}[1]{\boldsymbol{#1}}

\title{Solving the Brachistochrone Problem\\ by an Influence Diagram\thanks{This work 
was supported by the Czech Science Foundation (project 16-12010S).}}
\author{Ji\v{r}\'{\i} Vomlel}
\affil{Institute of Information Theory and Automation,\\
Czech Academy of Sciences, \\
Pod vod\'arenskou v\v{e}\v{z}\'{\i} 4, Prague 8, 182 08, Czechia\\
vomlel@utia.cas.cz,\\
\texttt{http://www.utia.cas.cz/vomlel/}}
\date{}

\begin{document}

\maketitle             

\begin{abstract}
Influence diagrams are a decision-theoretic extension of probabilistic graphical models.
In this paper we show how they can be used to solve the Brachistochrone problem. 
We present results of numerical experiments on this problem,
compare the solution provided by the influence diagram with the optimal solution.
The R code used for the experiments is presented in the Appendix.
\end{abstract}

\section{Introduction}

Formulated by Johan Bernoulli in 1696, the brachistochrone problem is: given two points find a curve connecting them such that 
a mass point moving along the curve under the gravity reaches the second
point in minimum time. See~\citep[Example 3.4.2]{bertsekas-2000} for a formulation of this problem 
as an optimal control problem.

\section{The ODE model}

The state variable is the vertical coordinate $y$. It is assumed to be a function
of the horizontal coordinate $x$. The control variable $u$ controls the derivative of $y$:
\begin{eqnarray*}
 \dfrac{d y(x)}{d x}  & = & u(x)
\end{eqnarray*}
The task is to find the control function $u(x)$ so that we get from a point $(0,0)$ to $(a,b)$, where 
$a > 0$ and $b < 0$.
This means that the boundary conditions are 
\begin{eqnarray*}
y(0) & = & 0 \\
y(a) & = & b \enspace .
\end{eqnarray*}
It is also assumed that the initial speed at the origin is zero.

Speed $v$ is defined by the law of energy conservation 
-- kinetic energy equals to the change of gravitational potential energy:
\begin{eqnarray}
\frac{1}{2} \cdot m \cdot (v)^2 & = & - m \cdot g \cdot y \\
v & = & \sqrt{-2 \cdot g \cdot y} \label{eq-v}\enspace .
\end{eqnarray}

For an infinitesimal segment of length $dx$ with an infinitesimal change $dy$ 
of the vertical position $y$ we can write
\begin{eqnarray}
v & = & \dfrac{ds}{dt} \ =  \ \sqrt{\dfrac{dy}{dt}^2+\dfrac{dx}{dt}^2} \ = \ 
\sqrt{\left(\dfrac{dy}{dx}\dfrac{dx}{dt}\right)^2+\dfrac{dx}{dt}^2} \ = \
\left(\sqrt{1+\dfrac{dy}{dx}^2}\right) \dfrac{dx}{dt}  \label{eq-v2} \enspace . \rule{3mm}{0mm}
\end{eqnarray}
By substituting~\eqref{eq-v} to~\eqref{eq-v2} we get
\begin{eqnarray}
dt & = & \dfrac{ds}{v} \ \ = \ \ \left(
\dfrac{1}{\sqrt{-2\cdot g \cdot y}}
\sqrt{1+\dfrac{dy}{dx}^2}\right) dx \enspace .
\end{eqnarray}

The solution of the brachistochrone problem is a function $y=f(x)$ that minimizes the total 
time $T$ necessary to get from the point $(0,0)$ to the point $(a,b)$ 
\begin{eqnarray}
T & = & \int_{0}^{a} \left(\dfrac{1}{\sqrt{-2\cdot g \cdot f(x)}} \sqrt{1+\dfrac{d f(x)}{dx}^2}\right) dx \enspace . \label{eq-time}
\end{eqnarray}

The solution of the brachistochrone problem is known -- it is a part of a cycloid, which
can be specified by parametric formulas:
\begin{eqnarray*}
x & = & \dfrac{K}{2} \left( \varphi + \sin \varphi \right) + L \\
y & = & - \dfrac{K}{2} \left( 1 - \cos \varphi \right) \enspace . 
\end{eqnarray*}
The constants $K,L$ are specified so that the cycloid goes trough points $(0,0)$ and $(a,b)$.

\section{Discretized version of the problem}

We discretize the problem:
\begin{itemize}
\item $n$ ... the number of discrete intervals, 
\item $\Delta x =\dfrac{a}{n}$ ... discretization step of the x-coordinate,
\item $i$ ... the index of the discrete interval,
\item $x_i$ ... x-coordinate $i \cdot \Delta x$, $i=0,1,\ldots,n$
\item $y_i$ ... y-coordinate at $x_i$,
\item $v_i$ ... speed at $x_i$,
\item $u_i$ ... control at coordinate $x_i$,
\item $t_{i}$ ... time to get from $x_{i-1}$ to $x_i$.
\end{itemize}

The state variable $y_i$ is transformed by the control variable $u_i$ as
\begin{eqnarray*}
y_{i+1} & = & y_{i} + u_{i}\enspace . 
\end{eqnarray*}
In each segment we will assume that the path is a line segment, i.e. for
$x \in [x_i,x_{i+1}]$ and for $y \in [y_i,y_{i+1}]$ it holds that
\begin{eqnarray}
y & = &  \dfrac{u_i}{\Delta x} \cdot x + y_i \label{eq-line} \enspace .
\end{eqnarray}
By substituting~\eqref{eq-line} to~\eqref{eq-time} and by solving the integral 
we get the formulas for the time spent at the segment $[x_i,x_{i+1}]$.
\begin{eqnarray}
t_{i+1} & = & \left\{\begin{array}{ll}
\dfrac{\Delta x}{\sqrt{-2 \cdot g \cdot y_i}}& \mbox{if $u_i=0$} \\[4mm]
-\sqrt{\dfrac{2}{g}} \cdot \left(\dfrac{(\Delta x)^2+u_i^2}{u_i} \right) \cdot
\left(\sqrt{-y_i} - \sqrt{-u_i-y_i}\right) & \mbox{otherwise.}
\end{array}\right.
\label{eq-ti}\enspace
\end{eqnarray}


The boundary conditions are
\begin{eqnarray*}
(x_0, y_0) & = & (0,0) \\
(x_n, y_n) & = & (a,b) \enspace .
\end{eqnarray*}

The goal is to find the control strategy $\bs{u}=(u_0,\ldots,u_{n-1})$, $u_i \in \mathbb{R}$, $i=0,1\ldots,n-1$ 
so that we get from the initial point $(x_0,y_0)$ to the terminal point $(x_n,y_n)$ in the shortest possible time
\begin{eqnarray*}
J(\bs{u}) & = & \sum_{i=1}^{n} t_i \enspace 
\end{eqnarray*}
and satisfy the state conditions (the gravitational potential energy corresponding to the value of $y$
cannot be more than it was at the initial point):
\begin{eqnarray*}
y_i & \leq & y_0 \ \ \mbox{for $i=1,\ldots,n$.} \enspace 
\end{eqnarray*}

\section{The influence diagram}

We will illustrate how an influence diagram can be used to find an arbitrary precise
solution of the problem. 
An influence diagram~\citep{howard-matheson-1981} is a Bayesian network
augmented with decision variables and utility functions. For details see, e.g., \cite{jensen-2001}.

\begin{figure}[htb]
\begin{center}
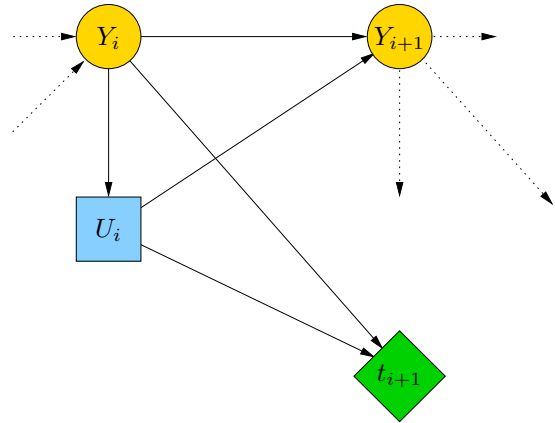
\caption{A Segment of the Influence Diagram for the Brachistochrone Problem\label{fig-id}}
\end{center}
\end{figure}

The structure of a segment of the influence diagram for the discrete version 
of the Brachistochrone Problem is presented in Figure~\ref{fig-id}.
The utility function for node $t_{i+1}$ is defined by formula~\eqref{eq-ti}.
The conditional probability $P(Y_{i+1}|U_i,Y_i)$ is deterministic and defined as:
\begin{eqnarray*}
P(Y_{i+1}=y_{i+1}|U_i=u_i,Y_i=y_i) & = & \left\{ \begin{array}{ll}
1 & \mbox{if $y_{i+1}=y_i+u_i$}\\
0 & \mbox{otherwise.}
\end{array}\right.
\end{eqnarray*}

\begin{figure}[htb]
\begin{center}
\includegraphics[width=0.9\textwidth]{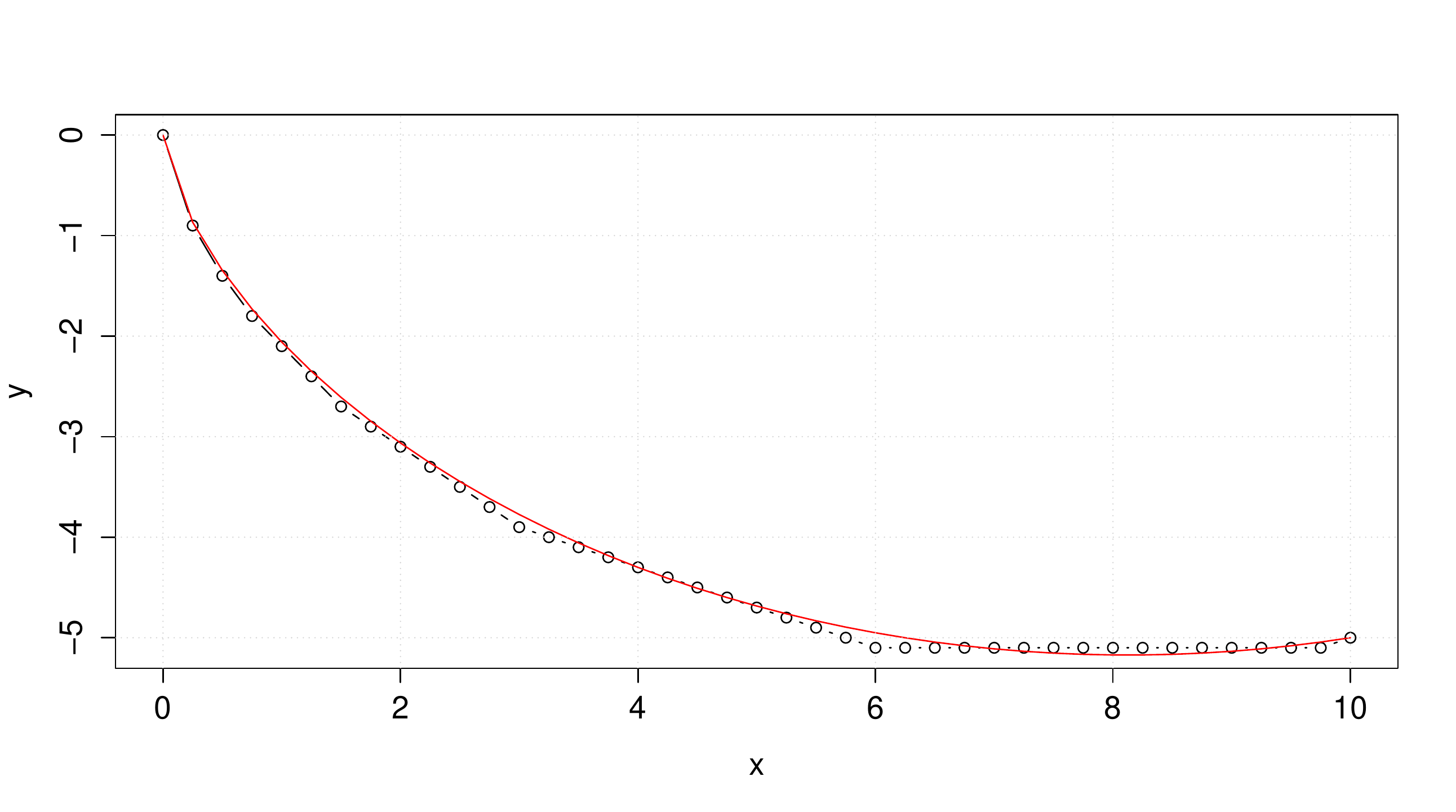}
\end{center}
\caption{Comparison of the optimal solution with the influence diagram solution.\label{fig-compare}}
\end{figure}

In Figure~\ref{fig-compare} we compare the optimal trajectory (full red line) with
the solution found by the influence diagram (circles connected by lines)
for $\Delta x = 0.25$, $\Delta y = 0.1$ and $(a,b)=(10,-5)$.
The difference between the optimal trajectory and the influence diagram solution 
can be reduced by reducing the discretization steps $\Delta x$ and $\Delta y$.
The experiments were performed using R~\citep{r-2014} 
-- we present the code in Appendix~\ref{R-code}.

\clearpage

\section{Conclusions}

We have shown how influence diagrams can be used to solve the Brachistochrone problem.
The numerical experiments reveal that the solution found by influence diagrams 
approximates well the optimal solution. In future we plan to apply influence 
diagrams to other trajectory optimization problems where the optimal solution is not known.
These problems are traditionally solved by methods of optimal control
theory but influence diagrams offer an alternative that can bring several
benefits over the traditional approaches.

\bibliographystyle{apalike}
\bibliography{bibliography}
 
\appendix

\section{The R code\label{R-code}}

\begin{verbatim}
n.x <- 41 # number of x values
n <- 101 # number of y values
a <- 10 # the x-coordinate of the goal state
b <- -5 # the y-coordinate of the goal state
delta.x <- a/(n.x-1) # the discretization step of x
delta.y <- 2*(-b)/(n-1) # the discretization step of y
g <- 9.81 # the gravitation constant
eps <- 10^-12 
  
# time spent at one segment of length delta.x assuming linear path 
time.step <- function(u,y){
  if ((y>0) || (y+u>0) || (((y==0)&(u==0)))){
    return(Inf) 
  }else{
    if (u==0){
      return(delta.x/sqrt(-2*g*y))
    }else{
      s <- sqrt(delta.x^2 + u^2)
      return(sqrt(2/g)*(s/u)*(sqrt(-y) - sqrt(-(u + y))))
    }
  }
}

address.y <- function(y){
  stopifnot(y <= 0)
  stopifnot(y >= 2*b)
  return(round(1+(y-2*b)/delta.y))
}
value.y <- function(address){
  return(2*b+(address-1)*delta.y)
}

address.u <- function(u){
  stopifnot(u <= -b)
  stopifnot(u >= b)
  return(round(1+(u-b)/delta.y))
}
value.u <- function(address){
  return(b+(address-1)*delta.y)
}

address.x <- function(x){
  stopifnot(x <= a+eps)
  stopifnot(x >= 0)
  return(round(1+(x/delta.x)))
}
value.x <- function(address){
  return((address-1)*delta.x)
}

is.addmissible <- function(y,k){
  if (k==(n.x-1)){
    return(abs(y-b)<eps)
  }else{
    return((y <= 0) & (y >= 2*b)) 
  }
}

find.best.policy <- function(y.start=0){
  policy <- array(0,dim=c(n.x-1,n))
  expected.utility <- rep(0,times=n)
  cat("\n")
  for (k in (n.x-1):1){
    x <- value.x(k)
    expected.utility.new <- rep(Inf,times=n)
    for (i in 1:n) {
      y <- value.y(i)  
      for(j in 1:n){    
        cat("\r k=",k," i=",i,"j=",j,"                                ")
        u <- value.u(j)
        y.next <- y+u
        # if y.next is within the admissible region 
        if (is.addmissible(y.next,k)){
          exp.util <- time.step(u=u,y=y) 
					            + expected.utility[address.y(y.next)]
          if (exp.util < expected.utility.new[address.y(y)]){
            expected.utility.new[address.y(y)] <- exp.util
            policy[address.x(x), address.y(y)] <- u
          }    
        }
      }  
    } 
    expected.utility <- expected.utility.new
  }
  return(list(policy=policy,
	            expected.utility=expected.utility[address.y(y.start)]))
}

# The construction of the state (vertical position y) profile.
# Note: since u and y have the same discretization step it is assured that
# by the application of u at state y we stay at the grid of y
construct.y.profile <- function(policy, y.start=0){
  x <- 0
  y <- y.start
  profile.y <- array(0,dim=c(n.x))
  profile.y[1] <- y
  for (i in 1:(n.x-1)){
    u <- policy[i,address.y(y)]
    y <- y+u
    profile.y[i+1] <- y  
  }
  return(profile.y)      
}

# The construction of the control profile.
# Note: since u and y have the same discretization step it is assured that
# by the application of u at state y we stay at the grid of y
construct.u.profile <- function(policy, y.start=0){
  x <- 0
  y <- y.start
  profile.u <- array(0,dim=c(n.x-1))
  for (i in 1:(n.x-1)){
    u <- policy[i,address.y(y)]
    y <- y+u  
    profile.u[i] <- u  
  }
  return(profile.u)      
}

evaluate.u.profile <- function(profile.u, y.start=0){
  val <- 0
  y <- y.start
  for(i in 1:length(profile.u)){
    u <- profile.u[i]
    val <- val + time.step(u,y)
    y <- y + u
  }
  return(val)
}

# Brachistochrone (the solution found by the Mathematica FindRoot function)
theta.max <- 3.50837
a.val <- 2.586
theta.val <- (0:100)*(theta.max/100)
brachistochrone.x <- a.val * (theta.val - sin(theta.val)) 
brachistochrone.y <- - a.val * (1 - cos(theta.val))

# The actual computations
res <- find.best.policy()
profile.y <- construct.y.profile(res$policy)

# Plot results
plot(x=(0:(n.x-1))*delta.x, y=profile.y, type="b", xlab="x", ylab="y")
lines(x=brachistochrone.x, y=brachistochrone.y, col="red")
grid()

profile.u <- construct.u.profile(res$policy)
plot(x=0:(n.x-2), y=profile.u, type="l", xlab="x", ylab="u")
grid()

# reconstruction of the optimal control profile 
# from the values found by the Mathematica FindRoot function
brachistochrone.x <- (0:40)*delta.x
# values found by the Mathematica FindRoot function
brachistochrone.y <- c(0,        -0.86755, -1.34679, -1.73084, -2.05946, 
                       -2.34941, -2.60981, -2.84634, -3.06283, -3.26204, 
                       -3.44602, -3.61637, -3.77433, -3.92094, -4.05702, 
                       -4.18327, -4.30028, -4.40854, -4.50848, -4.60047, 
                       -4.68481, -4.76179, -4.83164, -4.89456, -4.95074, 
                       -5.00031, -5.04342, -5.08018, -5.11067, -5.13497, 
                       -5.15315, -5.16523, -5.17125, -5.17123, -5.16517, 
                       -5.15304, -5.13483, -5.11049, -5.07995, -5.04316, -5.0)
brachistochrone.u <- -brachistochrone.y[-length(brachistochrone.y)]
                     +brachistochrone.y[-1]

# compute the total time for the path found by the influence diagram
evaluate.u.profile(profile.u)
# compute the total time for the optimal path at the same discrete scale
evaluate.u.profile(brachistochrone.u)
\end{verbatim}

\end{document}